\documentclass[11pt]{amsart}
\usepackage{amssymb,latexsym}

\newdimen\AAdi%
\newbox\AAbo%
\def\AAk#1#2{\s_etbox\AAbo=\hbox{#2}\AAdi=\wd\AAbo\kern#1\AAdi{}}%
\def\AAr#1#2#3{\s_etbox\AAbo=\hbox{#2}\AAdi=\ht\AAbo\raise#1\AAdi\hbox{#3}}%

\font\tenmsb=msbm10 at 12pt
\font\sevenmsb=msbm7 at 8pt
\font\fivemsb=msbm5 at 6pt
\newfam\msbfam
\textfont\msbfam=\tenmsb
\scriptfont\msbfam=\sevenmsb
\scriptscriptfont\msbfam=\fivemsb
\def\Bbb#1{{\tenmsb\fam\msbfam#1}}

\newcommand{\beq}{\begin{equation}}
\newcommand{\eeq}{\end{equation}}
\newcommand{\beqr}{\begin{eqnarray}}
\newcommand{\eeqr}{\end{eqnarray}}
\newcommand{\ba}{\begin{array}}
\newcommand{\ea}{\end{array}}

\begin{document}

\newtheorem{thm}{Theorem}
\newtheorem{lem}{Lemma}
\newtheorem{cor}{Corollary}
\newtheorem{rem}{Remark}
\newtheorem{pro}{Proposition}
\newtheorem{defi}{Definition}
\newtheorem{conj}[thm]{Conjecture}
\newcommand{\noi}{\noindent}
\newcommand{\dis}{\displaystyle}
\newcommand{\mint}{-\!\!\!\!\!\!\int}

\def \bx{\hspace{2.5mm}\rule{2.5mm}{2.5mm}} \def \vs{\vspace*{0.2cm}}
\def\hs{\hspace*{0.6cm}}
\def \ds{\displaystyle}
\def \p{\partial}
\def \O{\Omega}
\def \o{\omega}
\def \b{\beta}
\def \m{\mu}
\def \l{\lambda}
\def\L{\Lambda}
\def \ul{u_\lambda}
\def \D{\Delta}
\def \d{\delta}
\def \k{\kappa}
\def \s{\sigma}
\def \e{\varepsilon}
\def \a{\alpha}
\def \tf{\widetilde{f}}
\def\cqfd{%
\mbox{ }%
\nolinebreak%
\hfill%
\rule{2mm} {2mm}%
\medbreak%
\par%
}
\def \pr {\noindent {\it Proof.} }
\def \rmk {\noindent {\it Remark} }
\def \esp {\hspace{4mm}}
\def \dsp {\hspace{2mm}}
\def \ssp {\hspace{1mm}}

\def \u{u_+^{p^*}}
\def \ui{(u_+)^{p^*+1}}
\def \ul{(u^k)_+^{p^*}}
\def \energy{\int_{\R^n}\u }
\def \sk{\s_k}
\def \mo{\mu_k}
\def\cal{\mathcal}
\def \I{{\cal I}}
\def \J{{\cal J}}
\def \K{{\cal K}}
\def \OM{\overline{M}}
\def \L{{\cal L}_H}

\def\fk{{{\cal F}}_k}
\def\M1{{{\cal M}}_1}
\def\Fk{{\cal F}_k}
\def\Fl{{\cal F}_l}
\def\FF{\cal F}
\def\Gk{{\Gamma_k^+}}
\def\n{\nabla}
\def\uuu{{\n ^2 u+du\otimes du-\frac {|\n u|^2} 2 g_0+S_{g_0}}}
\def\uuug{{\n ^2 u+du\otimes du-\frac {|\n u|^2} 2 g+S_{g}}}
\def\sku{\sk\left(\uuu\right)}
\def\qed{\cqfd}
\def\vvv{{\frac{\n ^2 v} v -\frac {|\n v|^2} {2v^2} g_0+S_{g_0}}}
\def\vvs{{\frac{\n ^2 \widetilde v} {\widetilde v}
 -\frac {|\n \widetilde v|^2} {2\widetilde v^2} g_{S^n}+S_{g_{S^n}}}}
\def\skv{\sk\left(\vvv\right)}
\def\tr{\hbox{tr}}
\def\pO{\partial \Omega}
\def\dist{\hbox{dist}}
\def\RR{\Bbb R}\def\R{\Bbb R}
\def\C{\Bbb C}
\def\B{\Bbb B}
\def\N{\Bbb N}
\def\Q{\Bbb Q}
\def\Z{\Bbb Z}
\def\PP{\Bbb P}
\def\EE{\Bbb E}
\def\F{\Bbb F}
\def\G{\Bbb G}
\def\H{{\mathcal H}}
\def\SS{\Bbb S}\def\S{\Bbb S}

\def\lcf{{locally conformally flat} }

\def\circledwedge{\setbox0=\hbox{$\bigcirc$}\relax \mathbin {\hbox
to0pt{\raise.5pt\hbox to\wd0{\hfil $\wedge$\hfil}\hss}\box0 }}

\def\sss{\frac{\s_2}{\s_1}}

\def\dd{\delta_0}

\date{}
\title[A Hardy-Moser-Trudinger inequality]{A Hardy-Moser-Trudinger inequality}

\thanks{G.W. is partly supported by SFB/TR71 of DFG and 
D.Y. is supported by the French
ANR project referenced ANR-08-BLAN-0335-01. The work was started  during a visit of D.Y. at University of Freiburg and continued during a visit
of G.W. at Metz University. The authors would like to thank both
departments 
for their 
hospitality. G.W. would also like to thank the Department of Mathematics, University of Magdeburg for its hospitality.}
\author{Guofang Wang}
\address{ Albert-Ludwigs-Universit\"at Freiburg,
Mathematisches Institut,
Eckerstr. 1,
D-79104 Freiburg, Germany}
\email{guofang.wang@math.uni-freiburg.de}
\author{Dong Ye}
\address{LMAM, UMR 7122,
D\'epartement de Math\'ematiques,
Universit\'e Paul Verlaine de Metz, F-57045 Metz, France}
\email{dong.ye@univ-metz.fr}
\begin{abstract} In this paper we obtain an inequality on
 the unit disc $B$ in the plane, which improves the classical Moser-Trudinger inequality and
 the classical Hardy inequality at the same time. Namely, there exists a constant $C_0>0$ such that
\[
\int_B e^{\frac {4\pi u^2}{H(u)}} dx \le C_0 < \infty, \quad \forall\; u\in C^\infty_0(B),\]
where
 $$H(u) := \int_B |\n u|^2 dx - \int_B \frac {u^2}{(1-|x|^2)^2} dx.$$
  This inequality is a two dimensional analog of the
 Hardy-Sobolev-Maz'ya inequality in higher dimensions, which was recently intensively studied. We also prove that
 the supremum is achieved in a suitable function space, which is an analog of the celebrated
 result of Carleson-Chang for the Moser-Trudinger inequality.
\end{abstract}

\maketitle
\section{Introduction}
Let $B$ denote the standard unit disc in $\R^2$. The famous Moser-Trudinger inequality \cite{Moser, Trudinger}  
\begin{equation}
\label{MT}
\int_B e^{\frac {4\pi u^2}{\|\nabla u\|_2^2}} dx \le C < \infty, \quad \forall\; u \in H_0^1(B)\end{equation}
plays an important role in two dimensional analytic problems.
 This inequality is viewed as an analog of the Sobolev inequality for the higher dimensional cases. It is optimal in the  sense that
the constant $4\pi$ could not be replaced by any larger constant.
Its slighter weaker form 
\[\frac 12 \int_B |\nabla u|^2 dx - 8\pi \log \left(\int_B e^u dx\right) \ge -C > -\infty, \quad \forall\; u\in H^1_0(B)\]
has been  intensively used in the problem of prescribed Gaussian  curvature and recently in the mean field equation. Here again, the constant $8\pi$ is optimal to have a finite infimum.

\medskip
There is another important inequality in analysis, the Hardy inequality
\[H(u) := \int_B |\n u|^2 dx - \int_B \frac {u^2}{(1-|x|^2)^2} dx \ge 0,\quad \forall\; u\in H^1_0(B).\]
Here and after, $|\cdot|$ denotes always the Euclidean norm. This result is also optimal in the sense that for any $\lambda > 1$, 
$$\inf_{H_0^1(B)}\left[\int_B |\nabla u|^2 dx - \lambda\int_B \frac {u^2}{(1-|x|^2)^2} dx\right] = -\infty.$$
Indeed, this inequality holds also for higher dimensions. The Hardy inequality can be improved in the following way. Let $B^n$ denotes the unit ball in $\R^n$ with $n \geq 2$, it is known that there exists $C > 0$ (see \cite{BM} or Remark \ref{embed} below for $n = 2$) such that
\begin{align}
\label{refineH}
H(u) \geq C\int_{B^n} u^2 dx, \quad \forall \; u \in H_0^1(B^n).
\end{align}
Hence $\|u\|_\H := \sqrt{H(u)}$ defines a norm over $H_0^1(B^n)$ and the completion of $C_0^\infty(B^n)$ with respect to the norm $\|\cdot\|_\H$ is a Hilbert space,
which is denoted by $\H(B^n)$. Obviously $H_0^1(B) \subsetneq \H(B)$. For simplicity, we denote $\H(B)$ by $\H$ and $\|\cdot\|_\H$ by $\|\cdot\|$.

\medskip
In this paper, one of our main objectives is to improve the Moser-Trudinger inequality by combining the Hardy inequality. 
\begin{thm}\label{thm1}
There exists a constant $C_0>0$ such that 
\begin{equation}\label{HMT}
\int_B e^{\frac {4\pi u^2}{H(u)}} dx \le C_0 < \infty, \quad \forall\; u\in \H.\end{equation}
\end{thm}

\smallskip
A direct consequence is the following, slightly weaker, but applicable form.

\begin{cor} There exists a constant $C>0$ such that 
\begin{equation}\label{weakform}
\frac 12 \int_B |\n u|^2 dx - \frac 12 \int_B \frac {u^2}{(1-|x|^2)^2}dx - 8\pi \log \left(\int_B e^u dx\right) \ge -C > -\infty, \quad \forall\; u\in \H\backslash\{0\}.
\end{equation}
\end{cor}

\medskip
In the first glimpse these improved inequalities look too strong to be true. But if one compares to the recent work on the
Hardy-Sololev inequality in higher dimensional case, one
would speculate that this can be true. The Hardy-Sobolev inequality for higher dimension is also called the Hardy-Sobolev-Maz'ya inequality. Maz'ya  proved in \cite{M} (section 2.1.6, Corollary 3)
that there exists a constant $C > 0$ such that for any $u\in H^1_0(B^n)$ with $n>2$, 
 \begin{equation} \label{2a}
\int_{{B}^n} |\nabla u(x) |^2 d x
-\int_{{B}^n}\frac{u^2}{(1-|x|^2)^2} dx\ge C
 \left(\int_{{B}^n} |u(x)|^{\frac{2n}{n-2}} dx \right)^{\frac{n-2}{n}}.
\end{equation}
Let us denote $C_n$  the best constant such that (\ref{2a}) holds. In some recent works, the constant $C_n$ has been estimated (see \cite{TT} for $n > 3$ and \cite{BFL} for $n = 3$). Let $S_n$ be the best constant for the Sobolev embedding from $H^1(\R^n)$ into $L^{\frac{2n}{n - 2}}(\R^n)$.
\begin{itemize}
\item If $n>3$, then $C_n < S_n$ and the sharp constant $C_n$ is achieved in $\H(B^n)$. See \cite{TT}.
\item If $n=3$, then $C_3 = S_3$ is not achieved in $\H(B^n)$. See \cite{BFL}.
\end{itemize}
Theorem \ref{thm1} shows that the analog of the Hardy-Sobolev-Maz'ya inequality (\ref{2a}) holds on a two dimensional disc with the same best constant $4\pi$ as in the classical Moser-Trudinger inequality (\ref{MT}). We call (\ref{HMT}) or \eqref{weakform} a Hardy-Moser-Trudinger inequality.  Moreover, we obtain also a Carleson-Chang type result, i.e., the best constant is achieved.

\begin{thm}\label{thm2} There exists $u_0 \in \H$ such that $\|u_0\| = 1$ and
\[
\int_B e^{4\pi u_0^2} dx = \max_{u\in \H, \|u\| \leq 1}\int_B e^{4\pi u^2} dx=\max_{u\in \H\setminus\{0\}}\int_B e^{\frac {4\pi u^2}{H(u)}} dx.
\]
\end{thm}
Note that the supremum is not achieved in $H_0^1(B)$, see Remark \ref{rem3} below.

\medskip
We wonder if this kind of Hardy-Moser-Trudinger inequality  holds for more general domains $\O \subset \R^2$. For example, let $\O \subset \R^2$ be a regular, bounded and convex domain, then (see \cite{BM})
\begin{align*}
H_d(u) := \int_\O |\n u|^2 dx - \frac{1}{4}\int_B \frac {u^2}{d(x, \p\O)^2} dx > 0,\quad \forall\; u\in H^1_0(\O)\backslash\{0\}.
\end{align*}
We propose the following

\medskip
\noindent
{\bf Conjecture}: {\it There is a constant $C(\O)>0$ such that 
\begin{align*}
\int_\Omega e^{\frac {4\pi u^2}{H_d(u)}} dx \le C(\O) < \infty, \quad \forall\; u\in \H_d(\O)\setminus\{0\}.
\end{align*}}

\noindent
Here $\H_d(\O)$ denotes the completion of $C_0^\infty(\O)$ with the corresponding norm, defined by $\|u\|_{H_d}^2 = H_d(u)$. The conjecture is true when $\O = B$. This follows immediately from Theorem \ref{thm1}, since $H(u) \leq H_d(u)$ for any $u \in C^\infty_0(B)$.

\medskip
There is another  improved Moser-Trudinger inequality on the  disc in $\R^2$, which was recently
proved and studied in \cite{MS, AT}. 
\begin{align}
\label{MT_hyper}
\sup_{u\in H^1_0(B),\|\nabla u\|_2 \le 1 }\int_B \frac{e^{4\pi u^2} - 1}{(1-|x|^2)^2} dx < \infty.
\end{align}

For other generalizations  of the classical Moser-Trudinger inequality (\ref{MT}), see for instance
\cite{AD, FOR, LLY,LR}. See also \cite{DT, TZ} for a Moser-Trudinger inequality in K\"ahler geometry.
It would be interesting to know if there are such generalizations for 
our inequality (\ref{HMT}). A generalization to a higher dimensional ball like in \cite{Lin} and to higher order derivatives like in \cite{A, Beckner, F, LY} would be also very interesting. For the  related
higher order equations see also \cite{GGS}.

\medskip
The proof of the Hardy-Moser-Trudinger inequality \eqref{HMT} is  different from Moser's  proof for \eqref{MT}.
Here we use the blow-up analysis, which is an important tool in geometric analysis.
 A similar approach was also used in \cite{CC, DJLW, AD, Li, LR}
to establish 2-dimensional inequalities. The blow-up analysis for the elliptic equation related to the classical Moser-Trudinger inequality \eqref{MT} was initiated in \cite{S, AS}. It would be an interesting question if there
is a Moser-like proof for the Hardy-Moser-Trudinger inequality \eqref{HMT}.

\medskip
The paper is organized as follows: In Section 2, we reduce our problem  to radially nonincreasing functions
and study the property of such functions with bounded $H(u)$. We also study Green's function to the
operator
$-\Delta -\frac 1{(1-|x|^2)^2}$. In Section 3, we prove a subcritical Hardy-Moser-Trudinger inequality
with constant $(4\pi-\e)$ for any constant $\e \in (0, 4\pi)$. The subcritical Hardy-Moser-Trudinger inequality
(\ref{HMTe}) is achieved by a function $u_\e \in \H$. In Section 4, we analyze the convergence of the sequence $\{u_\e\}$ as $\e$ tends to $0$, and its blow-up behavior. Finally, Theorems \ref{thm1} and \ref{thm2} are proved respectively in Section 5 and 6, by contradiction arguments.

\medskip
In the following, $\|\cdot\|_p$ denotes the standard $L^p$ norm for $p \in [1, \infty]$ and $C$ denotes a generic positive constant, which could be changed from one line to another.

\section{Preliminaries}
First of all, we use the nonincreasing symmetrization with respect to the standard hyperbolic metric $dv_{\cal H} = \frac{dx}{(1 - |x|^2)^2}$ over $B$, which will enable us to reduce the problem on radially symmetric functions.

\medskip
For any $u \in H_0^1(B)$, let $u_*$ denote the associated radially nonincreasing rearrangement with respect to $dv_\H$, it is well known that $u_* \in H_0^1(B)$ and $\|\nabla u_*\|_2 \leq \|\nabla u\|_2$ (see \cite{Ba}). On the other hand,
$$\int_B \frac{u_*^2}{(1-|x|^2)^2} dx = \int_B  u_*^2 dv_{\cal H} = \int_B  u^2 dv_{\cal H} = \int_B \frac{u^2}{(1-|x|^2)^2} dx.$$ 
So $H(u) \leq 1$ implies that $H(u_*) \leq 1$. Furthermore,
\begin{align*}
\int_B e^{4\pi u_*^2} dx = \int_B e^{4\pi u_*^2}(1 - |x|^2)^2 dv_{\cal H} \geq \int_B e^{4\pi u^2}(1 - |x|^2)^2 dv_{\cal H} = \int_B e^{4\pi u^2} dx.
\end{align*}
This estimate comes from the Hardy-Littlewood inequality (see for example \cite{Brock}) by noticing that the rearrangement of $(1 - |x|^2)^2$ is just itself. Therefore we only need to consider nonincreasing, radially symmetric functions. Let
\begin{align*}
\Sigma :=\Big\{u \in C_0^\infty(B), u(x)=u(r) \mbox{ with $r=|x|$, $u'\leq 0$}\Big\}
\end{align*} and $\H_1$ be the closure of $\Sigma$ in $\H$. To prove Theorem \ref{thm1}, it remains to show that 
\begin{align}
\sup_{u \in \H_1, \|u\| \leq 1} \int_B e^{4\pi u^2} dx \leq C_0 < \infty.
\end{align}

\medskip
Let $u \in \Sigma$, $r=\varphi(t)=\tanh\left(\frac t 2\right)$ and 
\begin{align*}
\widetilde u (t) := u(r)=u\circ \varphi (t).
\end{align*}
Denote $B_r$ the disk of radius $r$ centered  at $0$ and $B^c_r=B\backslash B_r$ its compliment in $B$. Define also
\[H_\Omega(u) := \int_\Omega |\n u|^2 dx - \int_\Omega \frac {u^2}{(1-|x|^2)^2} dx , \quad \forall\; \Omega \subset B.\]
It is easy to see that
\[H_{B_r^c}(u) = 2\pi \int_{t = \varphi ^{-1}(r)}^\infty \left(\widetilde u'^2 - \frac{\widetilde u^2}{4}\right)\sinh(s) ds, \quad \forall\; r \in (0, 1).\]
Set now $v(s)= e^\frac{s}{2} \widetilde u(s)$. Note that  $v(s) = 0$ for large $s$. Integration by parts gives then
\begin{align}\label{add_x1}
\begin{split}
\frac{H_{B_r^c}(u)}{2\pi} & = \int_t^\infty e^{-s} \left(v'^2-v v'\right) \sinh(s) ds \\
&= \frac {v^2(t)} 2 e^{-t} \sinh(t) + \int_t^\infty e^{-2s} v^2 ds + \int_t^\infty e^{-s} v'^2\sinh(s) ds.
\end{split}
\end{align}
We obtain, by taking  $t$ tending to $0$, i.e. $r \rightarrow 0$,
\begin{align*}
\frac{H(u)}{2\pi} = \int_0^\infty e^{-2s} v^2 ds + \int_0^\infty e^{-s} v'^2\sinh(s) ds, \quad \forall\; u \in \Sigma.
\end{align*}
Consequently, we have
\begin{lem}
\label{Hrad}
$\H_1$ is embedded continuously in $H^1_{loc}(B) \cap C^{0, \frac{1}{2}}_{loc}(B\setminus\{0\})$.  Moreover, for any $p \geq 1$, $\H_1 \subset L^p(B)$ and these embeddings are compact.
\end{lem}
\pr Indeed, fix any $r \in (0, 1)$, there holds that for all $u \in \Sigma$,
\begin{align}
\label{newembed}
\begin{split}
\int_{B_r} u^2 dx + \int_{B_r} |\nabla u|^2 dx & = 2\pi\int_0^{\varphi^{-1}(r)} \left[\widetilde u^2 \frac{1}{4\cosh^4(\frac{s}{2})} + \widetilde u'^2\right] \sinh(s)ds\\
& \leq C\int_0^{\varphi^{-1}(r)} \big(\widetilde u^2 + \widetilde u'^2\big)\sinh(s) ds\\
& \leq C\int_0^{\varphi^{-1}(r)} \big(v^2 e^{-s} + v'^2e^{-s}\big)\sinh(s) ds\\
& \leq C_r\int_0^{\varphi^{-1}(r)} \big[v^2 e^{-2s} + v'^2e^{-s}\sinh(s)\big] ds\\
& \leq C_rH(u). 
\end{split}
\end{align}
Here the constants $C_r$ depend only on $r \in (0, 1)$. From above we have 
 $\H_1 \subset H^1_{loc}(B)$.
By the Sobolev embedding, we get
$$\H_1 \subset \cap_{p \geq 1}L^p_{loc}(B) \quad \mbox{and} \quad\H_1 \subset C^{0, \frac{1}{2}}_{loc}(B\setminus\{0\}).$$

\medskip
Furthermore, for any $r \in (0, 1)$, there exists $C_r > 0$ such that $u(r) \leq C_rH(u)$, $\forall\; u \in \H_1$. As $u \in \H_1$ is nonincreasing with respect to the radius, there holds that $\H_1$ is continuously embedded in $L^p(B)$ for any $p \geq 1$. Taking any bounded sequence $\{u_k\} \subset \H_1$, up to a subsequence we may assume that $u_k$ converges to $u$ weakly in $\H$ and  a.e. in $B$.  Finally we conclude by the following Lemma. \qed

\begin{lem}
\label{cvL1}
Let $\O \subset \R^n$ be of finite measure. If a sequence of measurable functions $w_k$ converges a.e. in $\O$ to $w$, and there exists $q > 1$ such that $\{w_k\}$ is bounded in $L^q(\O)$, then $w_k$ converges to $w$ in $L^1(\O)$.
\end{lem}

\begin{rem}
\label{embed}
Using the symmetrization argument, we can see that $\H$ is embedded continuously in $L^p(B)$ for any $p \in [1, \infty)$.
\end{rem}

 From (\ref{add_x1}) we have
\begin{align*}
\frac{H_{B_r^c}(u)}{2\pi} \geq \frac {v^2(t)} 2 e^{-t} \sinh(t) \geq 0, \quad \forall \; u \in \Sigma, \; r = \tanh\left(\frac{t}{2}\right),\; t > 0.
\end{align*}
Hence $H_{B_r}(u) \leq H(u)$ for any $u \in \Sigma$ and $r \in (0, 1)$. Pay attention that 
$H_{B^c_r}(u) \le H(u)$
 is in general not true. Since $\sinh(t) = \frac{2r}{1 - r^2}$, the above implies 
\begin{lem} \label{lem1} For  any  $u \in \H_1$
\begin{equation}\label{1}
u^2(r)\le \frac{1 - r^2}{2\pi r}H_{B_r^c}(u), \quad \forall\; r \in (0, 1].\end{equation}
\end{lem}

\medskip
Another crucial point in our approach is to handle the Hardy operator 
$$\L = -\Delta -\frac 1{(1-|x|^2)^2}$$ 
The problem is not trivial because we cannot apply directly the classical theory to $\H$ due to the potential which is singular on the boundary. Our idea is to separate the study into two parts, use the classical theory near the origin, then the $L^2$ theory in $\H$ for the exterior part. 
\begin{pro}
\label{solL1}
For any $f \in L^1(B)\cap L^2(B_{\frac{1}{8}}^c)$, there exists a unique 
\begin{align*}
v \in \H + W_0^{1, p}(B_\frac{1}{2})\;\; \mbox{with $p \in (1, 2)$}
\end{align*}
such that
\begin{align}
\label{fL1}
\L (v) = -\Delta v - \frac{v}{(1 - |x|^2)^2} = f \quad \mbox{in } {\cal D}'(B).
\end{align}
Moreover, we can decompose $v = v_1 + v_2$ where $v_1 \in \H$, $v_2 \in \cap_{p < 2}W^{1, p}_0(B_\frac{1}{2})$ and 
\begin{align*}
\|v_1\| + \|\nabla v_2\|_p \leq C_p\|f\|_1 + C\|f\|_{L^2(B_{\frac{1}{8}}^c)} \quad \forall \; p \in (1, 2).
\end{align*} 
\end{pro}
\begin{rem}
Of course, the decomposition $(v_1, v_2)$ is not unique, however the solution $v$ is uniquely determined. There exist some interesting works for $L^2$ theory with more general singular potentials, see for example \cite{DD}.
\end{rem}

\pr To simplify notations, define $\O_1 = B_{\frac{1}{8}}^c$, $\O_2 = B_{\frac{1}{2}}$ and 
\begin{align*}
a(x) = \frac{1}{(1 - |x|^2)^2}. \end{align*}
For the uniqueness of $v$, we need only to consider the case $f = 0$. Let $v = v_1 + v_2$ satisfy $\L(v) = 0$ in ${\cal D}'(B)$ with $v_1 \in \H$, $v_2 \in W^{1, p}_0(\O_2)$ and $p > 1$. We have
\begin{align*}
\left\langle v_1, \varphi\right\rangle_\H \leq C\|\varphi\| \leq C\|\nabla \varphi\|_{L^2(\O_2)}, \quad \forall\; \varphi \in C_0^\infty(\O_2).
\end{align*}
By equation, it implies that $-\Delta v_2 - a(x)v_2 \in H^{-1}(\O_2)$, the dual space of $H_0^1(\O_2)$. Given any $u \in H_0^1(\O_2)$, let $w(x) = u(\frac{x}{2}) \in H_0^1(B)$. Using the monotonicity of $a$, we have
\begin{align*}
0 \leq H(w)  \leq \int_{\O_2}|\nabla u|^2 dx - \int_B a(x)w^2 dx \leq \int_{\O_2}|\nabla u|^2 dx - 4\int_{\O_2}a(x)u^2 dx,
\end{align*}
and hence $H_{\O_2}(u) \geq \frac{3}{4}\|\nabla u\|_2^2$ for all $u \in H_0^1(\O_2)$. Therefore the operator $\L$ is coercive in $H_0^1(\O_2)$ and the classical regularity theory implies that $v_2$ belongs to $H_0^1(\O_2)\subset H_0^1(B)\subset \H$. Finally $v = v_1 + v_2 \in \H$ verifies $\left\langle v, \varphi\right\rangle_\H = 0$ for any $\varphi \in C_0^\infty(B)$. By density argument, $v = 0$, namely there is at most one solution.

\medskip
For the existence of solutions to \eqref{fL1} with $f \in L^1(B)\cap L^2(\O_1)$, consider first
\[\L(w) = f \;\; \mbox{in } \O_2, \quad w = 0 \;\; \mbox{on } \p \O_2.\]
 From the standard elliptic theory, there exists a unique solution $w \in \cap_{p < 2} W^{1, p}_0(\O_2)$ and $\|\nabla w\|_p \leq C_p\|f\|_1$. Choose $\Psi \in \Sigma$ a cut-off function such that $\Psi(r) = 1$ for $r \leq \frac{1}{8}$ and $\Psi(r) = 0$ for $r \geq \frac{1}{4}$. It is easy to check that  $$\L((1 - \Psi)w) = (1 - \Psi)f + 2\nabla w\nabla \Psi + w\Delta\Psi=:f_1\; \; \mbox{in } {\cal D}'(\O_2)$$ 
with $f_1 \in L^p(\O_2)$, $\forall\; p \in (1, 2)$. Thus  we get $(1 - \Psi)w \in W^{2, p}(\O_2)$, $\forall\; p \in (1, 2)$. In particular, applying the Sobolev embedding, $w$ extended by $0$ lies in $W^{1, q}(B_\frac{1}{4}^c)$ for all $q > 1$.

\medskip
Define now $\Psi_1(x) = \Psi(2x)$ and $h = (1 - \Psi_1)f + 2\nabla w\nabla \Psi_1 + w\Delta\Psi_1$. There holds $h \in L^2(B)$, hence $\|h\varphi\|_1 \leq \|h\|_2\|\varphi\|_2 \leq C\|h\|_2\|\varphi\|$ using \eqref{refineH} or Remark \ref{embed}. By the Riesz Theorem, there exists unique $v_1 \in \H$ such that 
\[\left\langle v_1, \varphi\right\rangle_\H = \int_B h\varphi dx, \quad \forall\; \varphi \in \H.\]
Clearly 
\begin{align*}
\|v_1\|\leq C\|h\|_2 \leq C\|f\|_{L^2(\O_1)} + C_p\|(1 -  \Psi)w\|_{W^{2, p}(\O_2)} \leq C_p\|f\|_1 + C\|f\|_{L^2(\O_1)}.
\end{align*}
Finally, let $v_2 = w\Psi_1$, we check readily that $v = v_1 + v_2$ is the desired solution. \qed

\smallskip
Using this result, we define a Green's function associated to the operator $\L$. 
\begin{pro}
\label{Gdef}
There exists a unique function $G_0 \in \H + W^{1, p}_0(B_\frac{1}{2})$ with $p \in [1, 2)$ such that
\begin{align}
\label{equa_G}
\L(G_0) = \delta_0 \quad \mbox{in } {\cal D}'(B)
\end{align}
where $\delta_0$ stands for the Dirac distribution at $0$. Moreover, $G_0$ is a radial function and there is a constant $C_G \in \R$ such that for any $\alpha \in (0, 1)$,
\begin{align}
\label{expG}
G_0(r) = -\frac{\ln r}{2\pi} + C_G + O\left(r^{1+\alpha}\right) \quad \mbox{as $r \to 0$},
\end{align}
\end{pro} 
\pr With a similar idea as above, let 
\begin{align*}
G_2(r) = -\frac{1}{2\pi}\Psi(r)\ln(r),\quad  F(r) =  -\frac{\ln r}{2\pi(1 - r^2)^2}\Psi - \frac{\Psi'}{\pi r} - \frac{\ln r}{2\pi}\D\Psi.
\end{align*}
Here $\Psi$ is the same cut-off function as in the previous proof. It is clear that $F \in L^2(B)$. Denote $G_1$ the unique solution in $\H$ such that 
\[\left\langle G_1, \varphi\right\rangle_\H = \int_B F\varphi dx, \quad \forall\; \varphi \in \H.\]
Clearly, $G_0 =  G_2  + G_1$ satisfies  equation \eqref{equa_G}. The uniqueness of $G_0$ is ensured by Proposition \ref{solL1}, which implies then $G_0$ is radial. For the expansion, since $F$ belongs to $L^p(B)$ for any $p > 1$, the standard elliptic theory yields that $G_1 \in W^{2, p}_{loc}(B) \subset C^{1, \alpha}_{loc}(B)$ for any $\alpha \in (0, 1)$. \qed

\begin{rem}
As $(1 - \Psi)\frac{\ln r}{2\pi} \in H_0^1(B)$, we have $G_0(r) = -\frac{\ln r}{2\pi} + \overline G(r)$ in $B$, with $\overline G \in \H$.
\end{rem}

\section{A weaker Hardy-Moser-Trudinger inequality}
In this section we prove a weaker form of the Hardy-Moser-Trudinger inequality, or its  subcritical version,
 which will be used in our proof of Theorem \ref{thm1}.
\begin{thm}\label{thm3}
For any constant $\e\in (0, 4\pi)$, it holds
\begin{align}
\label{HMTe}
 \sup_{u \in \H_1, \|u\| \leq 1} \int_B e^{(4\pi - \e)u^2} dx < \infty.
\end{align}
and the supremum is achieved by some $u_\e \in \H_1$. \end{thm}

\medskip
Define 
\begin{align*}
A_u(r) = \frac{1}{\pi r^2}\int_{B_r}  \frac {u^2}{(1-|x|^2)^2} dx.
\end{align*}
We estimate $A_u(r)$ for general $u \in \H_1$.
\begin{lem}
\label{lemA}
Let $u \in \H_1$ and $r \in (0, 1)$, we have 
\begin{align*}
\pi \left(\frac 1 2-r^2\right)A_u(r) \leq H(u) + \frac {\pi u(r)^2}{1-r^2}.
\end{align*}
\end{lem}
\pr By an elementary inequality
\[(u-b)^2+b^2 \ge \frac{u^2}{2}, \quad \forall\; u, b\in \R,\]
we deduce, for any $r \in (0, 1)$,
\begin{align*}
\frac 1{r^2} \int_{B_r} \frac {(u-b)^2}{(1-|x|^2)^2} dx & \geq
\frac 1{2r^2} \int_{B_{r}} \frac {u^2}{(1-|x|^2)^2} dx -\frac {1}{r^2} \int_{B_{r}} \frac {b^2}{(1-|x|^2)^2} dx \\
&= \frac {\pi} 2A_u(r)-\frac {\pi b^2}{1-r^2}.
\end{align*}
Applying the above to $w(x) = u(rx) - u(r) \in H_0^1(B)$ and $b = u(r)$, we get
\[
\int_B\frac{w^2}{(1-|x|^2)^2} dx\ge \int_{B_{r}}\frac {(u- u(r))^2}{r^2(1-|x|^2)^2} dx\ge \frac \pi 2 A_u(r)-\frac{\pi u(r)^2}{1-r^2}.\]
It follows that, together with the Hardy inequality,
\begin{align*}
0 \le H(w)  = \int_{B_r} |\n u|^2 dx - \int_B \frac{w^2}{(1-|x|^2)^2} dx & \le \int_{B_r} |\n u|^2 dx - \frac{\pi}2 A_u(r)+\frac{\pi u(r)^2}{1-r^2}\\
&=  H_{B_r} (u) + \pi r^2A_u(r) - \frac{\pi} 2A_u(r) +\frac{\pi u(r)^2}{1-r^2}\\
& \leq H(u)+\frac{\pi u(r)^2}{1-r^2} -\pi \left(\frac{1}{2}-r^2\right)A_u(r).
\end{align*}
which is just the conclusion. \qed

\medskip
\noindent{\it Proof of Theorem \ref{thm3}.}
Let $u \in \Sigma$ with $H(u) \leq 1$ satisfy
$$ \int_B e^{(4\pi - \e)u^2} dx \geq \sup_{w \in \H_1, \|u\| \leq 1} \int_B e^{(4\pi - \e)w^2} dx - 1.$$
Assume now $u(0) > 1$, otherwise $\|e^{(4\pi - \e)u^2}\|_1 \leq \pi e^{(4\pi - \e)}$ and \eqref{HMTe} holds true. Define \[r_1 = \inf\{r > 0 \;| \; u(r) \leq 1\} > 0.\]
Since $u(r_1) = 1$, we have a upper bound of $r_1$, $ r_1 \leq \frac{1}{2\pi}$ by \eqref{1}. From  Lemma \ref{lemA} and the fact $u(r) \geq 1$ for $r \leq r_1$, there exists $C > 0$  independent of $u$ such that $A_u(r) \leq Cu(r)^2$, $\forall\; r \in (0, r_1]$. Using again the estimate \eqref{1}, we have, for any $r \leq r_1$, 
\begin{align}
\label{estA2}
\begin{split}
\int_{B_r}|\n u|^2 dx = H_{B_r}(u)+ \pi r^2A_u(r) & \le 1- H_ {B_r^c}(u)  +  C\pi r^2 u(r)^2\\
& \leq 1 -  H_ {B_r^c}(u) + \frac{C}{2}r H_ {B_r^c}(u).
\end{split}
\end{align}
Therefore, there exists $r_2 \in (0, r_1]$ small enough, independent of $u$,
 such that $$\|\nabla u\|_{L^2(B_{r_2})} \leq 1.$$Thanks to the
 Moser-Trudinger inequality (\ref{MT}), 
\begin{align*}
\int_{B_{r_2}} e^{4\pi [u(r) - u(r_2)]^2} dx = \int_B e^{4\pi [u(r) - u(r_2)]_+^2} dx \leq C_{MT} < \infty.
\end{align*}
In thefollowing, we denote $f_+=\max\{f, 0\}$.
Moreover, by  Lemma \ref{lem1} and  $H(u) \leq 1$,  it holds $u(r_2) \leq C$ for some constant $C>0$ independent of $u$.
Hence we have that for any $r \leq r_2$,
\begin{align*}
(4\pi -\e)u(r)^2 \leq 4\pi \big[u(r) - u(r_2)\big]^2 + 8\pi u(r)u(r_2) - \e u(r)^2 \leq 4\pi \big[u(r) - u(r_2)\big]^2 + C_\e.
\end{align*}
Here $C_\e$ depends only on $\e$. This yields
\begin{align*}
\int_B e^{(4\pi - \e)u^2} dx & = \int_{B_{r_2}} e^{(4\pi - \e)u^2} dx + \int_{B_{r_2}^c} e^{(4\pi - \e)u^2} dx \\ &\leq \int_{B_{r_2}} e^{4\pi[u - u(r_2)]^2 + C_\e} dx + \pi(1 - r_2^2) e^{(4\pi - \e) u(r_2)^2}\\
& \leq e^{C_\e}C_{MT} + \pi e^{4\pi C} < \infty.
\end{align*}
The proof of inequality \eqref{HMTe} is completed. 

\medskip
Now we show the achievement of the supremum. Fix $\e > 0$, consider a maximizing sequence $u_j \in \H_1$ for \eqref{HMTe} with $\|u_j\|\leq 1$. Recall that $\|\cdot\|$ is the norm in $\H$.
Up to a subsequence, we can assume that $u_j$ converges weakly in $\H$ to $u_\e \in \H_1$, so 
$\|u_\e\| \leq 1$. By Lemma \ref{Hrad}, we can assume also that $u_j$ converges to $u_\e$ a.e. in $B$.

\medskip
Using \eqref{HMTe} with $\frac{\e}{2}$, we see that $e^{(4\pi - \e)u_j^2}$ is bounded in $L^q(B)$ for some $q > 1$. 
Lemma \ref{cvL1} implies then $e^{(4\pi - \e)u_j^2}$ converges in $L^1(B)$, that is,
\[\int_B e^{(4\pi - \e)u_\e^2} dx = \lim_{j\to\infty}\int_B e^{(4\pi - \e)u_j^2} dx,\]
and hence the supermum of \eqref{HMTe} is attained by $u_\e$. Clearly we must have $\|u_\e\|= 1$. \qed

\section{Blow-up analysis}
For any $\e \in (0, 4\pi)$, let $u_\e$ realize the maximum obtained by Theorem \ref{thm3}.
In this section we consider the convergence of the sequence $\{u_\e\}$ when $\e$ goes to zero, or more precisely, we try to understand the behavior of $\{u_\e\}$ if $\|u_\e\|_\infty$ explodes.

\medskip
Suppose that $\|u_ \e\|_\infty = u_\e(0)$ does not go to infinity as $\e$ tends to $0$, then there exists $\e_j \to 0$ such that $\|u_{\e_j}\|_\infty \leq C$. It is easy to see that in this case, up to a subsequence
 $u_{\e_j}$ converges weakly to $u_0 \in \H_1$ in $\H$ and a.e. in $B$, so $\|u_0\| \leq 1$ and $u_0 \in L^\infty(B)$. Let $w \in \H$, $\|w\| \leq 1$, 
\begin{align*}
 \int_B e^{(4\pi - \e_j)w^2} dx \leq  \int_B e^{(4\pi - \e_j)u_{\e_j}^2} dx, \quad \mbox{for any }\; j \in \N.
\end{align*}
Applying respectively monotone and dominated convergence, it holds
\begin{align*}
 \int_B e^{4\pi w^2} dx  = \lim_{j\to \infty} \int_B e^{(4\pi - \e_j)w^2} dx \leq  \lim_{j\to \infty}\int_B e^{(4\pi - \e_j)u_{\e_j}^2} dx = \int_B e^{4\pi u_0^2} dx < \infty.
\end{align*}
In other words, $u_0$ realizes the finite maximum of the Hardy-Moser-Trudinger functional, therefore both Theorem \ref{thm2} and \ref{thm1} are proved in this case.

\medskip
In the following, we will suppose the contrary, i.e. $\lim_{\e \to 0} \|u_ \e\|_\infty = \infty$ and perform  a
blow-up analysis as in \cite{AD, LR}.
Since $u_\e \in \H_1$ is a maximizer, there exists $\l_\e > 0$ such that
\begin{align}
\label{equa_ue}
\L(u_\e) = -\D u_\e - \frac{u_\e}{(1-|x|^2)^2} = \l_\e u_\e e^{(4\pi - \e)u_\e^2} \quad \mbox{in } {\cal D}'(B).
\end{align}
Using Lemma \ref{Hrad} and \eqref{HMTe}, $\D u_\e \in L^q_{loc}(B)$ for some $q \in (1, 2)$, we get $u_\e \in W^{2, q}_{loc}(B)$ using the standard regularity theory. By the Sobolev embedding in dimension two, $u_\e$ is continuous in $B$, and hence $u_\e \in C(\overline B)$. Here the continuity up to $\p B$ follows from $u_\e \in \H_1$.

\medskip
Indeed, $u_\e$ is the so called $\H-$solution of \eqref{equa_ue} over $B$, in the spirit of D\'avila and Dupaigne \cite{DD}. Using $u_\e$ as a test function, we have
\begin{align}
\label{ue1}
\l_\e\int_B u_\e^2 e^{(4\pi - \e)u_\e^2} dx = \|u_\e\|^2 = 1.
\end{align}
\begin{rem}\label{rem3} Notice that $u_\e$ does not belong to $H_0^1(B)$.  This is due to Theorem III in \cite{BM}, because $u_\e \in C(\overline B)$ and
\begin{align*}
a(x) - \frac{1}{4d(x, \p B)^2} = \frac{1}{(1 - r^2)^2} - \frac{1}{4(1 - r)^2} = \frac{1}{1 - r} \times\frac{(3+r)}{4(1 + r)^2} = O\left(d(x, \p B)^{-1}\right).
\end{align*}
\end{rem}

Suppose that Theorem \ref{thm1} does not hold true, then
\begin{align}
\label{noHMT}
\lim_{\e\to 0} \int_B e^{(4\pi - \e)u_\e^2} dx = \infty = \lim_{\e\to 0} \|u_\e\|_\infty.
\end{align}
Since $\|u_\e\| = 1$, there exist weakly convergent subsequences in $\H$. Note that from now on, for simplicity, we do not distinguish always between convergence and subconvergence. Assume $u_\e \rightarrow u_0 \in \H_1$ weakly in $\H$. 
\begin{lem}
\label{u0=0}
We have $u_0 \equiv 0$.\end{lem}
\pr Suppose the contrary, then there is $r_0 \in (0, \frac{1}{2})$ such that $u_0(r_0) > 0$. 

\smallskip
By Lemma \ref{Hrad}, $u_\e$ tends to $u_0$ in $C_{loc}(B \setminus\{0\})$ (since the embedding of $C^{0, \frac{1}{2}}$ into $C^0$ is compact), hence $u_\e(r_0) \geq \delta > 0$ for $\e$ small enough. Using Lemma \ref{lemA}, when $\e$ is small enough, $A_{u_\e}(r) \leq Cu_\e(r)^2$ for any $r \leq r_0$, because $u_\e(r) \geq u_\e(r_0) \geq \delta$ and $H(u_\e) = 1$. Hence
\begin{align*}
\int_{B_{r}} |\nabla u_\e|^2dx = 1 - H_{B_{r}^c}(u_\e) + \pi {r}^2A_{u_\e}({r}) \leq 1 - \frac{2\pi r}{1 - r^2}u_\e(r)^2 + Cr^2u_\e(r)^2,\quad \forall\; r \leq r_0.
\end{align*}
There exists then $r_1 \in (0, r_0)$ and $\eta > 0$ such that for $\e$ small, $\|\nabla u_\e\|_{L^2(B_{r_1})} \leq 1 - \eta < 1$. By the Moser-Trudinger inequality (\ref{MT}),
\begin{align*}
\int_{B_{r_1}} e^{\frac{4\pi }{1 - \eta}(u_\e - b_\e)^2} dx \leq C_{MT}, \quad \mbox{where } b_\e = u_\e(r_1).
\end{align*}
Similarly as in the proof for \eqref{HMTe}, using $\lim_{\e \to 0} b_\e = u_0(r_1) < \infty$, we can conclude that $\|e^{4\pi u_\e^2}\|_1 \leq C < \infty$ for small enough $\e$, but this contradicts \eqref{noHMT}. 
\qed

\medskip
Applying Lemma \ref{u0=0} and Lemma \ref{Hrad}, $u_\e \in \H_1$ converges uniformly to $0$ in $\overline{B_r^c}$ for $r > 0$. Thus we will concentrate now our attention on the behavior of $u_\e$ near the origin. Define
\begin{align}\label{add1}
M_\e = u_ \e(0) = \max u_\e\quad \mbox{and}\quad r_\e^2 = \frac{e^{(\e - 4\pi)M_\e^2}}{\l_\e M_\e^2}.
\end{align}
Using \eqref{HMTe},
\begin{align*}
\l_\e^{-1} = \int_B u_\e^2 e^{(4\pi - \e)u_\e^2} dx \leq  \int_{B_\frac{1}{2}} u_\e^2 e^{(4\pi - \e)u_\e^2} dx + C & \leq M_\e^2 e^{(2\pi - \e)M_\e^2} \int_B e^{2\pi u_\e^2} dx + C\\
& \leq M_\e^2 e^{(2\pi - \e)M_\e^2} \sup_{\|u\| \leq 1}\int_B e^{2\pi u^2} dx + C\\
& \leq C M_\e^2 e^{(2\pi - \e)M_\e^2} + C.
\end{align*}
Consequently (recall that $\lim_{\e \to 0} M_\e = \infty$)
\begin{align*}
r_\e^2 M_\e^2 =  \frac{e^{(\e - 4\pi)M_\e^2}}{\l_\e} \leq CM_\e^2 e^{-2\pi M_\e^2},
\end{align*} 
hence $\lim_{\e\to 0} r_\e M_\e = 0$ and $\lim_{\e\to 0} r_\e = 0$. Define $v_\e(x) = u_\e(r_\e x)$ and $\xi_\e(x) = M_\e\big[v_\e(x) - M_\e\big]$, a direct calculation leads to
\begin{align}
\label{xi_e}
-\D \xi_\e =  \frac{v_\e}{M_\e} e^{(4\pi - \e)(v_\e^2 - M_\e^2)} + \frac{r_\e^2M_\e^2}{(1 - r_\e^2|x|^2)^2} \frac{v_\e}{M_\e}
\quad \mbox{in } {\cal D}'(B_{r_\e^{-1}}).
\end{align} 
For any $R > 0$, $-\D \xi_\e = O(1)$ in $B_R$ for small $\e$, since $0 \leq v_\e \leq M_\e$. By $\xi_\e(0) = 0$, the standard elliptic estimate implies that $\xi_\e$ converges in $C^1_{loc}(\R^2)$ to $\xi$. Therefore
\begin{align}
\label{est_ve}
v_\e - M_\e = \frac{\xi_\e}{M_\e} \rightarrow 0,\quad \frac{v_\e}{M_\e} \rightarrow 1 \quad\mbox{and}\quad  v_\e^2 - M_\e^2 = 2\xi_\e + \frac{\xi_\e^2}{M_\e^2} \rightarrow 2\xi \quad \mbox{in $C^1_{loc}(\R^2)$}.
\end{align}
Passing the limit $\e \to 0$ in \eqref{xi_e}, the equation verified by $\xi \in C^1(\R^2)$ is
\begin{align}\label{CL}
-\D \xi = e^{8\pi \xi} \quad \mbox{in } {\cal D}'({\R^2}).
\end{align}
Combining with the facts $\xi(0) = 0$, $\xi$ is radially symmetric and nonincreasing with respect to $r$, $\xi$ is uniquely determined,
\begin{align}
\label{equa_xi}
\xi(x) = -\frac{1}{4\pi}\ln(1 + \pi r^2), \quad \int_{\R^2}  e^{8\pi \xi} dx = 1.
\end{align}
Note that all solutions of (\ref{CL}) with $e^{8\pi \xi} \in L^1(\R^2)$ were classified in \cite{CL}.

\medskip
The above analysis is for understanding the behavior of the sequence $\{u_\e\}$ near the blow-up point $0$, more precisely in $B_{r_\e R}$ for any large, but fixed $R>0$. Let $L > 1$ and $R>0$ large. We divide the disc $B$ into three parts: the interior part $B_{r_\e R}$, the outer part $$\{Lu_\e \leq M_\e\}:=\left\{x\in B\,|\, u_\e(x)\le \frac{M_\e}{L}\right\}$$ 
and the neck region $$\{Lu_\e \geq M_\e\}\setminus B_{r_\e R} :=
\left\{x\in B\backslash B_{r_\e R}\,|\, u_\e(x)\ge \frac{M_\e}{L}\right\}.$$
To analyze $\{u_\e\}$ in the outer part and the neck region, let us denote
$u_{\e, L} = \min(u_\e, \frac{M_\e}{L})$.  
We have
\begin{lem}
\label{lemL}
For any $L > 1$,  $\limsup_{\e\to 0} H(u_{\e, L}) \leq L^{-1}$.
\end{lem}
\pr Consider $\zeta_{\e, L} = u_\e - u_{\e, L} = \left(u_\e - \frac{M_\e}{L}\right)_+$. Fix $R > 0$. Using $\zeta_{\e, L}$ as a test function to equation \eqref{equa_ue},  we have
\begin{align}
\label{add2}
\begin{split}
\int_B |\nabla\zeta_{\e, L}|^2 dx - \int_B \frac{u_\e\zeta_{\e, L}}{(1 - |x|^2)^{2}} dx&  = \l_\e\int_B\zeta_{\e, L}u_\e e^{(4\pi -  \e)u_\e^2}dx\\ & \geq \l_\e\int_{B_{r_\e R}} \zeta_{\e, L}u_\e e^{(4\pi -  \e)u_\e^2}dx\\
& = \int_{B_R}\left(\frac{v_\e}{M_\e} - \frac{1}{L}\right)_+ \frac{v_\e}{M_\e}e^{(4\pi -  \e)(v_\e^2 - M_\e^2)}dx\\
& \rightarrow \left(1 - \frac{1}{L}\right)\int_{B_R} e^{8\pi\xi} dx,
\end{split}
\end{align}
when $\e \to 0$. Recall that $r_\e$ is defined by (\ref{add1}), the convergence is ensured by \eqref{est_ve}. Moreover, one can check easily that
\begin{align*}
H(u_{\e, L}) = H(u_\e) - \int_B |\nabla\zeta_{\e, L}|^2 dx + \int_B \frac{u_\e\zeta_{\e, L}}{(1 - |x|^2)^2} dx,
\end{align*}
which, taking $R\to \infty$, together with \eqref{add2} and 
 \eqref{equa_xi}, completes the proof.
\qed

Using our subcritical inequality \eqref{HMTe} to functions $\frac{Lu_{\e, L}}{2}$, we get
\begin{align}
\label{ueL}
\int_B e^{\pi L^2u_{\e, L}^2} dx \leq C <  \infty, \quad \mbox{for $\e$ small enough}.
\end{align}
Furthermore,
\begin{lem}
\label{2limits} We have
\begin{align}
\label{limit}
\lim_{\e\to 0} \l_\e M_\e^2=0
\end{align}
and
\begin{align}
\label{Dirac}
\lim_{\e\to 0} \l_\e M_\e\int_B u_\e e^{(4\pi - \e)u_\e^2}dx = 1.
\end{align}
\end{lem}
\pr Let us estimate firstly $\|e^{(4\pi - \e)u_\e^2} \|_1$. Fix $L > 2$.
\begin{align*}
I_{\e} := \int_{\{Lu_\e \leq M_\e\}} e^{(4\pi - \e)u_\e^2} dx \leq \int_B e^{(4\pi - \e)u_{\e, L}^2} dx \rightarrow \pi, \quad \mbox{as $\e \to 0$}.
\end{align*}
The convergence follows from  the facts that $u_{\e, L} \to 0$ a.e. in $B$, the estimate \eqref{ueL} and Lemma
 \ref{cvL1}. Using once again the uniform convergence of $u_\e$  to zero in $B_r^c$ for any $r \in (0, 1)$, 
\begin{align*}
\int_{\{Lu_\e \leq M_\e\}} e^{(4\pi - \e)u_\e^2} dx \geq \int_{B_r^c} e^{(4\pi - \e)u_\e^2} dx \rightarrow \pi(1 - r^2), \quad \mbox{as $\e \to 0$}.
\end{align*}
Taking $r \to 0$, we can claim that $\lim_{\e \to 0} I_\e = \pi$. On the other hand, 
\begin{align*}
J_\e := \int_{\{Lu_\e \geq M_\e\}} e^{(4\pi - \e)u_\e^2} dx & \leq \frac{L^2}{\l_\e M_\e^2} \int_{\{Lu_\e \geq M_\e\}} \l_\e u_\e^2 e^{(4\pi - \e)u_\e^2} dx\\ & \leq  \frac{L^2}{\l_\e M_\e^2} \int_B \l_\e u_\e^2 e^{(4\pi - \e)u_\e^2} dx\\
& = \frac{L^2}{\l_\e M_\e^2}. 
\end{align*}
Finally,
\begin{align}\label{add3}
\infty = \lim_{\e \to 0} \int_B e^{(4\pi - \e)u_\e^2} dx & = \lim_{\e\to 0} \left(I_\e + J_\e\right) \leq \pi + \limsup_{\e\to 0}  \frac{L^2}{\l_\e M_\e^2},
\end{align}
which implies  $\liminf_{\e\to 0} \l_\e M^2_\e = 0$. This argument  is in fast valid for any subsequence. Hence we have \eqref{limit}.

\medskip
To prove \eqref{Dirac}, we estimate the integral over three parts separately. First
\begin{align*}
\l_\e M_\e \int_{\{Lu_\e \leq M_\e\}} u_\e e^{(4\pi - \e)u_\e^2} dx \leq \l_\e M_\e^2I_\e \rightarrow 0, \quad \mbox{as $\e \to 0$}.
\end{align*}
Moreover, for any $R > 0$,
\begin{align*}
\l_\e M_\e \int_{B_{r_\e R}} u_\e e^{(4\pi - \e)u_\e^2} dx = \int_{B_R}\frac{v_\e}{M_\e}e^{(4\pi -  \e)(v_\e^2  - M_\e^2)}dx \rightarrow \int_{B_R} e^{8\pi\xi} dx,
\end{align*}
and also
\begin{align*}
\l_\e M_\e \int_{\{Lu_\e \geq M_\e\}\setminus B_{r_\e R}} u_\e e^{(4\pi - \e)u_\e^2} dx & \leq L\int_{\{Lu_\e \geq M_\e\}\setminus B_{r_\e R}} \l_\e u_\e^2 e^{(4\pi - \e)u_\e^2} dx\\
& \leq L\int_{B\setminus B_{r_\e R}}\l_\e u_\e^2e^{(4\pi - \e)u_\e^2}dx\\
& = L - L\int_{B_{r_\e R}}\l_\e u_\e^2e^{(4\pi - \e)u_\e^2}dx\\
& \rightarrow L\left(1 - \int_{B_R} e^{8\pi\xi} dx\right).
\end{align*}
The proof of \eqref{Dirac} is completed by tending $R$ to $\infty$. \qed

\medskip
Let $g_\e = M_\e u_\e$, so $g_\e$ satisfies the equation 
\begin{align}\label{add4}
\L(g_\e) = \l_ \e g_\e e^{(4\pi - \e)u_\e^2} \quad \mbox{in } {\cal D}'(B).
\end{align}
\eqref{Dirac} and its proof shows that $\l_ \e g_\e e^{(4\pi - \e)u_\e^2}$ converges to the Dirac operator $\delta_0$ in ${\cal D}'(B)$, this suggests then $g_\e$ should tend to the corresponding Green's function $G_0$, which is confirmed as follows.
\begin{pro}
\label{neck}
The family $g_\e$ converges to $G_0$ in $W^{1, p}_{loc}(B)$ weakly if $p \in (1, 2)$, strongly in $L^q(B)$ for all $q \geq 1$ and also in $C(\overline{B_r^c})$, $\forall\; r \in (0, 1)$. Here $G_0$ is defined by Proposition
\ref{Gdef}.
\end{pro}
\pr Using Proposition \ref{solL1} on \eqref{add4}, as $g_\e \in \H$, there exist $k_\e$ and $h_\e$ such that $g_\e = h_\e + k_\e$ with $h_\e \in \H$, $k_\e \in \cap_{p < 2}W_0^{1, p}(B_\frac{1}{2})$  satisfying for any $p \in (1, 2)$,
\begin{align*}
\|h_\e\| + \|\nabla k_\e\|_p \leq C_p\|\l_\e g_\e e^{(4\pi - \e)u_\e^2}\|_1 + C\|\l_\e g_\e e^{(4\pi - \e)u_\e^2}\|_{L^2(\O_1)}\quad \mbox{where } \O_1 = B_\frac{1}{8}^c.
\end{align*}
As $u_\e e^{(4\pi - \e)u_\e^2}$ tends to zero uniformly in $\O_1$, so $\|h_\e\| + \|\nabla k_\e\|_p$ are uniformly bounded for $\e$ small, thus we get (up to subsequence) that $h_\e$ converges weakly to $h_0$ in $\H$ and $k_\e$ converges weakly to $k_0$ in $W_0^{1, p}(B_\frac{1}{2})$ for $p \in (1, 2)$. 

\medskip
Let $g_0 = h_0 + k_0 \in \H + W^{1, p}_0(B_\frac{1}{2})$, as $\l_ \e g_\e e^{(4\pi - \e)u_\e^2} \rightarrow \delta_0$, we have $\L(g_0) = \delta_0$ in ${\cal D}'(B)$. Applying Proposition \ref{Gdef}, $g_0 = G_0$ and we obtain all the claimed convergence of $g_\e$ to $G_0$. The convergence holds for the whole family $g_\e$ since the analysis is valid for any subsequence. \qed

\section{Proof of Theorem \ref{thm1}}
\noindent{\it Proof of Theorem \ref{thm1}.}
Suppose that Theorem \ref{thm1} does not hold. Let $\rho \in (0, 1)$ be small, which will be determined later. Thanks to Proposition \ref{neck}, From the previous section, recalling that $a(x) = (1 - |x|^2)^{-2}$, we have
\begin{align}
\label{cvge}
\lim_{\e\to 0} M_\e u_\e(\rho) = G_0(\rho), \quad \lim_{\e\to 0}\int_{B_{\rho}}a(x)M_\e^2u_\e^2 dx = \int_{B_{\rho}} a(x)G_0^2 dx =: J_1(\rho).
\end{align}
By equation \eqref{add4},  
\begin{align*}
H_{B_\rho^c}(g_\e) = H(g_ \e) - H_{B_\rho}(g_\e) = \l_\e \int_{B_\rho^c}g_\e e^{(4\pi - \e)u_\e^2} dx -\int_{\p B_\rho} \frac{\p g_\e}{\p \nu}g_\e d\sigma.
\end{align*}
Clearly as $\e \to 0$, the first term goes to $0$. By \eqref{Dirac} and \eqref{cvge},
\begin{align*}
 -\int_{\p B_\rho}  \frac{\p g_\e}{\p \nu}g_\e d\sigma = - g_\e(\rho)\int_{B_\rho}\D g_\e dx & =
 g_\e(\rho)\left(\int_{B_\rho} a(x)g_\e dx + \l_\e \int_{B_\rho}g_\e e^{(4\pi - \e)u_\e^2} dx\right)\\
& \rightarrow G_0(\rho)\left(\int_{B_\rho} a(x)G_0 dx + 1\right)=:J_2(\rho),
\end{align*}
Finally we have
\begin{align}
\label{estf1}
\int_{B_\rho} |\nabla u_\e|^2 dx = 1 - H_{B_\rho^c}(u_\e) + \int_{B_{\rho}}a(x) u_\e^2 dx = 1 - \frac{1}{M_\e^2}\Big[J_2(\rho) - J_1(\rho) + o_\e(1) \Big]
\end{align}
where $\lim_{\e\to 0} o_\e(1) = 0$, for any fixed $\rho > 0$. 

\medskip
Furthermore, using the expansion of $G_0$, direct calculations show that 
\begin{align*}
J_2(\rho) \sim -\frac{\ln\rho}{2\pi} \quad \mbox{and}\quad  J_1(\rho) \sim \frac{\rho^2(\ln\rho)^2}{4\pi}  \quad \mbox{as $\rho \to 0$}.
\end{align*}
So there is $\rho > 0$ small enough such that $J_2(\rho) - J_1(\rho)  > 0$. Fix such a $\rho$, the equality \eqref{estf1} implies that 
\begin{align*}
 \int_{B_\rho} |\nabla u_\e|^2 dx <  1 \quad \mbox{for $\e$ small enough}. 
\end{align*}
Applying the classical Moser-Trudinger inequality \eqref{MT}  to  $[u_\e - u_\e(\rho)]_+ \in H_0^1(B)$, we get
\[\int_{B_\rho} e^{4\pi[u_\e - u_\e(\rho)]^2} dx = \int_{B} e^{4\pi[u_\e - u_\e(\rho)]_+^2} dx\leq C_{MT}.\]
On the other hand, we have
\[u_\e^2(r) = [u_\e(r) - u_\e(\rho)]^2 + 2u_\e(r) u_\e(\rho) - u_\e^2(r) \leq [u_\e(r) - u_\e(\rho)]^2 + 2M_\e u_\e(\rho).\]
Therefore, when $\e$ tends to $0$, 
\begin{align*}
\int_B e^{4\pi u_\e^2} dx = \int_{B_\rho} e^{4\pi u_\e^2} dx + \int_{B_\rho^c} e^{4\pi u_\e^2} dx & \leq \int_{B_\rho} e^{4\pi [u_\e - u_\e(\rho)]^2 + 8\pi M_\e u_\e(\rho)} dx + \pi e^{4\pi u_\e(\rho)^2}\\
& \leq e^{8\pi g_\e(\rho)}C_{MT} + \pi e^{4\pi u_\e(\rho)^2}\\
& \rightarrow e^{8\pi G_0(\rho)}C_{MT} + \pi.
\end{align*} 
This contradicts obviously the hypothesis \eqref{noHMT}, and  hence the Hardy-Moser-Trudinger inequality must hold true. \qed

\section{Proof of Theorem \ref{thm2}}
\noindent{\it Proof of Thereom \ref{thm2}.}
Let $u_\e$ be the maximizer given by Theorem \ref{thm3}. We will proceed still by a contradiction argument
and borrow some ideas from \cite{CC, AD, LR}. Using the discussion at the end of Section 3, suppose that Theorem \ref{thm2} was not valid, then $\lim_{\e\to 0} \|u_\e\|_\infty = \infty$. 

\medskip
Define 
\begin{align*}
T_\e = \int_B e^{(4\pi - \e)u_\e^2}dx = \max_{u \in \H, \|u\|\leq 1} \int_B e^{(4\pi - \e)u^2}dx, \quad \forall\; \e \in (0, 4\pi).
\end{align*}
Clearly, $T_\e$ is increasing and 
\begin{align*}
\lim_{\e \to 0}T_\e = \sup_{u \in \H, \|u\|\leq 1} \int_B e^{4\pi u^2}dx := T_0.
\end{align*}
By Theorem \ref{thm1}, $T_0<\infty$. Readily $T_0>\pi$. 

\medskip
All arguments and properties obtained for $u_\e$ in the previous
section are true, except two points. One is the proof of the fact that  the weak limit $u_0$ is $0$ and 
another is the property
\eqref{limit}. 

\medskip
For the former point one can argue as follows. Fixing $\rho \in (0, 1)$, we see that for any $L > 1$, $u_{\e, L} = u_\e$ in $\overline{B_\rho^c}$ for $\e$ small enough, because $u_\e$ is uniformly bounded in $\overline{B_\rho^c}$ by Lemma \ref{Hrad} and $\lim_{\e\to 0}M_\e = \infty$. It follows, together with Lemma \ref{lemL}, 
\begin{align*}
\limsup_{\e \to 0} \|u_\e\|_{L^\infty(B_\rho^c)} = \limsup_{\e \to 0} \|u_{\e, L}\|_{L^\infty(B_\rho^c)} \leq C_\rho \limsup_{\e \to 0} H(u_{\e, L}) \leq \frac{C_\rho}{L}.
\end{align*}
Tend $L$ to $\infty$, we have $u_0 = 0$ in $B_\rho^c$. Since $\rho > 0$ is arbitrary, thus  $u_0 = 0$.

\medskip
The latter is no longer true. In fact, now we have
\begin{lem}
\label{lMe2} $\lim_{\e\to 0}\l_\e M_\e^2 = (T_0 - \pi)^{-1}$.
\end{lem}
\pr 
By the same argument as in the proof of \eqref{add3}, we get
\begin{align*}
T_0 = \lim_{\e \to 0}T_\e \leq \pi + \limsup_{\e\to 0}\frac{L^2}{\l_\e M_\e^2}, \quad \forall \; L > 2,
\end{align*}
which implies that $\limsup_{\e\to 0}\l_\e M_\e^2 < \infty$ since $T_0 > \pi$. Hence $\lim_{\e\to 0}\l_\e M_\e = 0$. Let $p_1 = \frac{L^2}{4} > 1$ (as $L > 2$). The estimate \eqref{ueL} and $\|u_\e\|_q \to 0$ for any $q \geq 1$ imply that as $\e \to 0$,
\begin{align*}
\int_{\{Lu_\e \leq M_\e\}} u_\e e^{(4\pi - \e)u_\e^2} dx \leq \|e^{(4\pi - \e)u_{\e, L}^2}\|_{p_1}\|u_\e\|_{q_1} \to 0 \quad \mbox{where } \frac{1}{p_1} + \frac{1}{q_1} = 1.
\end{align*}
The same argument shows that
\begin{align*}
\l_\e M_\e\int_{\{Lu_\e \leq M_\e\}} u_\e e^{(4\pi - \e)u_\e^2} dx \to 0, \quad \mbox{as $\e \to 0$}
\end{align*}
and \eqref{Dirac} holds true, so Proposition \ref{neck} remains true. Let $\O_{\rho, r} = B_\rho\setminus \overline{B_{r}}$ with $0 < r < \rho < 1$, as $\L(G_0) = 0$ in $\O_{\rho, r} \subset \R^2$, the Pohozaev identity shows that (recall that $a(x) = (1 - |x|^2)^{-2}$)
\begin{align*}
\int_{\O_{\rho, r}} \frac{{\rm div}\big[a(x)x]}{2}G_0^2(x) dx -\pi\Big[s^2G_0'^2(s) + a(s)s^2G_0^2(s)\Big]_r^\rho = 0.
\end{align*}
Using the expansion \eqref{expG} and tending $r \to 0$, we have
\begin{align*}
\int_{B_\rho} \frac{{\rm div}\big[a(x)x]}{2}G_0^2(x) dx -\pi\rho^2G_0'^2(\rho) -\pi a(\rho)\rho^2G_0^2(\rho) = -\frac{1}{4\pi}, \quad \forall\; \rho \in (0, 1).
\end{align*}
Similarly, applying the Pohozaev identity to $\L(u_\e) = \l_\e u_\e e^{(4\pi - \e)u_\e^2}$ in $B_\rho$ and multiplying by $M_\e^2$, we obtain that for any $\rho \in (0, 1)$,
\begin{align*}
& \int_{B_\rho} \frac{{\rm div}\big[a(x)x]}{2}g_\e^2(x) dx -\pi\rho^2g_\e'^2(\rho) -\pi a(\rho)\rho^2g_\e^2(\rho)\\ = & \; \l_\e M_\e^2 \left[\pi \rho^2\frac{e^{(4\pi - \e)u_\e(\rho)^2}}{4\pi - \e} - \int_{B_\rho} \frac{e^{(4\pi - \e)u_\e^2}}{4\pi - \e} dx\right].
\end{align*}
Finally, as $g_\e$ converges to $G_0$ in $C^1_{loc}(B\setminus\{0\})$ and $L^2(B)$ by standard elliptic theory and Proposition \ref{neck}, we obtain, for any $\rho \in (0, 1)$ (as $\limsup_{\e\to 0}\l_\e M_\e^2 < \infty$),
\begin{align*}
\l_\e M_\e^2\left[\int_{B_\rho} e^{(4\pi - \e)u_\e^2} dx - \pi\rho^2\right] \to 1,\quad \mbox{as $\e \to 0$}.
\end{align*}
Taking $\rho\to 1$ and using the uniform convergence of $u_\e$ to $0$ in $\overline{B_r^c}$ for $r > 0$,
 we are done. \qed

To get a contradiction, we proceed as in \cite{CC, LR}.  We first claim a Carleson-Chang type result.
\begin{lem}
\label{cc1}
If $\lim_{\e \to 0} \|u_\e\|_\infty = \infty$, then $T_0 \leq \pi(1 + e^{1 + 4\pi C_G})$ where $C_G$ is given by \eqref{expG}.
\end{lem}
\pr Fix $L > 2$ and let $R > 0$. Using \eqref{ueL} and Lemma \ref{cvL1}, we have the estimate for the exterior region as follows.
\begin{align*}
\int_{\{Lu_\e \leq M_\e\}} e^{(4\pi - \e)u_\e^2} dx \leq \int_B e^{(4\pi - \e)u_{\e, L}^2} dx \to \pi, \quad \mbox{as $\e \to 0$}.
\end{align*}
On the neck region, there holds
\begin{align*}
\int_{\{Lu_\e \geq M_\e\}\setminus B_{r_\e R}} e^{(4\pi - \e)u_\e^2} dx & \leq \frac{L^2}{\l_\e M_\e^2}\int_{\{Lu_\e \geq M_\e\}\setminus B_{r_\e R}} \l_\e u_\e^2 e^{(4\pi - \e) u_{\e, L}^2} dx \\
& \leq \frac{L^2}{\l_\e M_\e^2}\int_{B\setminus B_{r_\e R}}\l_\e u_\e^2 e^{(4\pi - \e)u_{\e, L}^2} dx\\
& = \frac{L^2}{\l_\e M_\e^2}\left(1 - \int_{B_{r_\e R}}\l_\e u_\e^2 e^{(4\pi - \e)u_{\e, L}^2} dx\right)\\
& \rightarrow L^2(T_0 - \pi)\left(1 - \int_{B_R} e^{8\pi\xi} dx\right).
\end{align*}
For  the integral over the interior region $B_{r_\e R}$, fix a small constant $\rho \in (0, 1)$. By \eqref{estf1},
we know that
\begin{align*}
\int_{B_\rho} |\nabla u_\e|^2 dx = 1 - \frac{E_\rho + o_\e(1)}{M_\e^2}
\end{align*}
where $\lim_{\e\to 0} o_\e(1) = 0$ and
\begin{align}
\label{E}
E_\rho :=J_2(\rho)-J_1(\rho)= G_0(\rho) + G_0(\rho)\int_{B_\rho} a(x)G_0 dx- \int_{B_\rho} a(x)G_0^2 dx.
\end{align}
Let $\ell_\e(x) = \|\nabla u_\e\|_{L^2(B_\rho)}^{-1}\big[u_\e(x) - u_\e(\rho)\big]_+$. Clearly, $\ell_\e \in H_0^1(B_\rho)$, $\|\nabla \ell_\e\|_2 = 1$ and $\ell_\e$ converges weakly to $0$ in ${\cal D}'(B_\rho)$. Using a  result of Carleson-Chang \cite{CC}, there holds
\begin{align}
\label{cc}
\limsup_{\e \to 0} \int_{B_\rho} \left(e^{4\pi \ell_\e^2} - 1\right) dx \leq \pi e\rho^2.
\end{align}
Moreover, $M_\e^{-1}\ell_\e  \to 1$ uniformly in $B_{r_\e R}$ since $ M_\e^{-1} u_\e \to 1$ uniformly in $B_R$ by \eqref{est_ve}. Therefore 
\begin{align*}
u_\e^2 (x) & = \big[\ell_\e(x) + u_\e(\rho)\big]^2\|\nabla u_\e\|_{L^2(B_\rho)}^2\\ & = \big[\ell_\e(x) + M_\e^{-1}G_0(\rho) + o_\e(M_\e^{-1})\big]^2\times\big[1 - M_\e^{-2}E_\rho + o_\e(M_\e^{-2})\big]\\
& = \ell_\e^2 (x)+ 2\ell_\e (x)M_\e^{-1} G_0(\rho) - \ell_\e^2(x)M_\e^{-2}E_\rho + o_\e(1)\\
& = \ell_\e^2(x) + 2G_0(\rho) - E_\rho + o_\e(1),
\end{align*}
where $o_\e(1)$ tends to $0$ uniformly in $B_{r_\e R}$. It implies that, together with \eqref{cc}, 
\begin{align*}
\limsup_{\e \to 0} \int_{B_{r_\e R}} e^{(4\pi - \e)u_\e^2} dx &\leq \limsup_{\e \to 0} \int_{B_{r_\e R}} \left(e^{4\pi u_\e^2} - 1\right)dx\\
& \leq e^{8\pi G_0(\rho) - 4\pi E_\rho}\limsup_{\e \to 0} \int_{B_{r_\e R}} \left(e^{4\pi\ell_\e^2} - 1\right)dx\\
& \leq e^{8\pi G_0(\rho) - 4\pi E_\rho}\limsup_{\e \to 0} \int_{B_\rho} \left(e^{4\pi\ell_\e^2} - 1\right) dx\\
& \leq \pi \rho^2e^{1+8\pi G_0(\rho) - 4\pi E_\rho}.
\end{align*}
Combining the three parts of estimation and let $R$ go to $\infty$, we conclude
\begin{align*}
T_0 = \lim_{\e\to 0}\int_B e^{(4\pi - \e)u_\e^2} dx \leq \pi + \pi \rho^2e^{1+8\pi G_0(\rho) - 4\pi E_\rho}, \quad \mbox{for any small $\rho > 0$}.
\end{align*}
Using the expansion \eqref{expG},
\[\frac 1{2\pi} \ln \rho + 2G_0(\rho)-E_\rho \to C_G, \quad \hbox{as }\; \rho \to 0.
\]
Hence it follows $T_0  \leq \pi(1 + e^{1 + 4\pi C_G})$. \qed

\medskip
We complete the proof of Theorem \ref{thm2} with the following lower bound estimate, which contradicts Lemma \ref{cc1}.
\begin{lem}
\label{lowerbT0}
There holds $T_0 > \pi(1 + e^{1 + 4\pi C_G})$.
\end{lem}
\pr The proof is a direct verification by choosing suitable test functions as in \cite{CC}.
Thanks to the blow-up analysis, we will consider a family $f_\e$ such that $f_\e$ looks like $M_\e^{-1}G_0$ outside a very small region of $0$ and $M_\e^{-1}\xi(r^{-1}_\e x) + M_\e$ near the origin where $\xi$ is given by \eqref{equa_xi}. For $\e > 0$ small, define
$$
f_\e(r) = \left\{ \begin{array}{ll}\beta_\e + \ds \vspace*{0,1cm} \frac{\xi(\e^{-1}r) + \gamma_\e}{\beta_\e}& \mbox{ if $r \leq \e R_\e$}\\
\ds \frac{G_0(r)}{\beta_\e}& \mbox{ if $\e R_\e \leq r \leq 1$}
\end{array}\right.
\quad \mbox{with } \; R_\e = -\ln\e.$$
Here  $\beta_\e$ and $\gamma_\e$ are constants to be chosen later. First, choose $\gamma_\e$ such that 
\begin{align*}
\frac{G_0(\e R_\e)}{\beta_\e} = \beta_\e + \frac{\xi(R_\e) + \gamma_\e}{\beta_\e},
\end{align*}
which makes functions $f_\e$ continuous.
Using the expansion of $G_0$, as $\e \to 0$, we have
\begin{align}
\label{fit}
\begin{split}
4\pi\left(\beta_\e^2 + \gamma_\e\right) & = -2\ln(\e R_\e) + 4\pi C_G + \ln(1 + \pi R_\e^2) + o\left(\e R_\e\right)\\
& =  -2\ln\e + 4\pi C_G + \ln\pi + O\left(R_\e^{-2}\right).
\end{split}
\end{align}
Clearly, $f_\e \in \H$. Now we estimate $\|f_\e\|$. Let $0 < r < \rho < 1$, by the equation of $G_0$,
\begin{align*}
H_{B_{\rho}\setminus B_{r}}(G_0) = \int_{\p (B_{\rho}\setminus B_{r})} G_0\frac{\p G_0}{\p \nu} d\sigma = -2\pi rG_0(r)G_0'(r) + 2\pi \rho G_0(\rho)G_0'(\rho) \leq -2\pi rG_0(r)G_0'(r),
\end{align*}
since $G_0$ is decreasing by the comparison principle. Taking $\rho \to 1$, we get
\begin{align*}
H_{B_{\e R_\e}^c}(f_\e) = \frac{H_{B_{\e R_\e}^c}(G_0)}{\beta_\e^2} & \leq -\frac{2\pi \e R_\e}{\beta_\e^2} G_0(\e R_\e)G_0'(\e R_\e)\\
& = -\frac{2\pi \e R_\e}{\beta_\e^2} \int_{\p B_{\e R_\e}} \frac{\p G_0}{\p \nu}d\sigma\\
& = \frac{G_0(\e R_\e)}{\beta_\e^2} \left(1 + \int_{B_{\e R_\e}} a(x) G_0 dx\right)\\
& = \frac{1}{4\pi\beta_\e^2}\big[-2\ln(\e R_\e) + 4\pi C_G + o\left(\e R_\e\right)\big].
\end{align*}
On the other hand, we have
\begin{align*}
\int_{B_{\e R_\e}} |\nabla f_\e|^2 dx = \frac{1}{\beta_ \e^2}\int_{B_{\e R_\e}} |\nabla \xi(\e^{-1}x)|^2 dx & = \frac{1}{\beta_ \e^2}\int_{B_R}|\nabla \xi|^2 dx\\
& = \frac{1}{4\pi\beta_ \e^2}\left[\ln(1 + \pi R_\e^2) - 1 + \frac{1}{1 + \pi R_\e^2}\right],
\end{align*}
and hence
\begin{align}
\label{testf}
H(f_\e) \leq  H_{B_{\e R_\e}^c}(f_\e)  + \int_{B_{\e R_\e}} |\nabla f_\e|^2 dx \leq \frac{1}{4\pi\beta_\e^2}\left[-2\ln\e + 4\pi C_G - 1 + \ln\pi + O\left(R_\e^{-2}\right)\right].
\end{align}
We choose $\beta_\e > 0$ such that $H(f_\e) = 1$, the estimate \eqref{testf} leads to (recall that $R_\e = -\ln\e$)
\begin{align*}
4\pi\beta_\e^2 \leq -2\ln \e + 4\pi C_G - 1 + \ln\pi + O\left(R_\e^{-2}\right) = O\left(|\ln \e|\right), \quad \mbox{as  $\e \to 0$}.
\end{align*}
From \eqref{fit}, it follows
\begin{align}
\label{bg}
4\pi\gamma_\e \geq 1 + O\left(R_\e^{-2}\right), \quad  \beta_\e = O\left(|\ln\e|^\frac{1}{2}\right)  \quad \mbox{as  $\e \to 0$}.
\end{align}

\smallskip
Now we estimate $\|e^{4\pi f_\e^2}\|_1$. Using $e^t \geq 1 + t$ in $\R$ and \eqref{bg}, we have
\begin{align*}
\int_{B_{\e R_\e}^c} e^{4\pi f_\e^2} dx  \geq \pi - \pi(\e R_\e)^2 + \frac{4\pi}{\beta_\e^2}\int_{B_{\e R_\e}^c} G_0^2 dx & = \pi + \frac{4\pi}{\beta_\e^2}\left[\int_B G_0^2 dx +  o\left(\e R_\e\right) - \frac{\e^2 R_\e^2\beta_\e^2}{4}\right]\\
& = \pi + \frac{4\pi}{\beta_\e^2}\left[\int_B G_0^2 dx + o_\e(1)\right].
\end{align*}
Moreover, in $B_{\e R_\e}$, by \eqref{fit} and  \eqref{bg} we have
\begin{align*}
4\pi f_\e^2(r) = 4\pi\left(\beta_\e + \frac{\xi(\e^{-1}r) + \gamma_\e}{\beta_\e}\right)^2 & \geq 4\pi\beta_\e^2 + 8\pi \gamma_\e + 8\pi\xi(\e^{-1}r)\\
& = 4\pi(\beta_\e^2 + \gamma_\e) + 4\pi\gamma_\e + 8\pi\xi(\e^{-1}r)\\
& \geq -2\ln\e + 4\pi C_G + \ln\pi + 1 + 8\pi\xi(\e^{-1}r) + O\left(R_\e^{-2}\right).
\end{align*}
The estimate is uniform in $B_{\e R_\e}$. Consequently
\begin{align*}
\int_{B_{\e R_\e}} e^{4\pi f_\e^2} dx  & \geq e^{-2\ln\e + 4\pi C_G + \ln\pi + 1 + O\left(R_\e^{-2}\right)} \int_{B_{\e R_\e}} e^{8\pi\xi(\e^{-1}r)} dx\\
& = \pi e^{4\pi C_G + 1 + O\left(R_\e^{-2}\right)}\int_{B_{R_\e}} e^{8\pi\xi} dx\\
& = \pi e^{4\pi C_G + 1}\left[1 + O\left(R_\e^{-2}\right)\right]
\end{align*}
where \eqref{equa_xi} was used in the last equality. Finally as $R_\e^{-2}\beta_\e^2 = o_\e(1)$, it holds
\begin{align*}
\int_{B} e^{4\pi f_\e^2} dx & \geq \pi + \frac{4\pi}{\beta_\e^2}\left[\int_B G_0^2 dx + o_\e(1)\right] + \pi e^{4\pi C_G + 1}\left[1 + O\left(R_\e^{-2}\right)\right]\\
&  = \pi + \pi e^{4\pi C_G + 1} +  \frac{4\pi}{\beta_\e^2}\left[\int_B G_0^2 dx + o_\e(1)\right].
\end{align*}
By choosing a small $\e>0$, we conclude readily $T_0 \geq \|e^{4\pi f_\e^2}\|_1 >  \pi + \pi e^{4\pi C_G + 1}$
Hence we finish the proof of Lemma \ref{lowerbT0}, and hence the proof of Theorem \ref{thm2}.\qed

\begin{rem}
 Like $u_\e$, the maximizer $u_0\in \H_1$ given by Theorem \ref{thm2} cannot belong to $H_0^1(B)$.
\end{rem}

\end{document}